# Адаптивная сегментация границ при квазистатическом анализе многопроводных линий передачи методом моментов

А. Е. Максимов, С. П. Куксенко

**Аннотация.** Описаны особенности квазистатического анализа многопроводных линий передачи методом моментов. Выполнено сравнение эффективности методов построения адаптивных сеток на примере вычисления ёмкостных матриц сложных многопроводных линий передачи. Выявлены оптимальные параметры методов, дающие уменьшение вычислительных затрат. Показано, что их использование даёт существенный прирост производительности при многовариантном анализе многопроводных линий передачи относительной густых равномерных сеток.

**Ключевые слова:** адаптивная сегментация границ, квазистатический анализ, метод моментов, многопроводные линии передачи, ёмкостная матрица.

## Введение

Как известно, основными элементами РЭС являются СВЧ-устройства, и прежде всего, линии передачи и резонаторы, поэтому требуется их качественное проектирование. В частности, при анализе линий передачи, в случае, когда поперечные размеры рассматриваемой структуры малы по сравнению с длиной распространяющейся электромагнитной волны, применяется квазистатический подход [1]. Он широко используется при проектировании многопроводных линий передачи (МПЛП). Особенностью такого проектирования является учёт распределенных первичных погонных параметров между всеми проводниками. По мере усложнения конструкций линий аналитические методы становятся всё менее применимы и поэтому прибегают к использованию численных методов, например методу моментов.

Для численного анализа какой-либо физической задачи необходимо построить её математическую модель, учитывающую существенные для данной задачи особенности реального объекта. Процесс построения математической модели для анализа электромагнитных задач формально можно представить в виде нескольких взаимосвязанных этапов, наиболее вычислительно затратными из которых являются построение сетки, формирование системы линейных алгебраических уравнений (СЛАУ) и её решение [2]. Так, выбор метода построения сетки оказывает влияние на свойства результирующей СЛАУ, что определяет выбор способа вычисления её элементов и метода её решения и тем самым временные затраты на её формирование и решение [3]. При необходимости проведения многовариантного анализа или оптимизации параметров исследуемой МПЛП эти три этапа повторяются многократно, что существенно усложняет процесс её оптимального проектирования из-за существенного роста вычислительных затрат. Цель данной работы – исследовать методы адаптивной сегментации границ МПЛП при квазистатическом анализе методом моментов.

**1. Математическая модель для вычисления ёмкостной матрицы МПЛП**

При квазистатическом подходе электрические характеристики, меняющиеся вдоль отрезков МПЛП длиной $dx$, характеризуются матрицами первичных погонных параметров **R** (Ом/м), **L** (Гн/м), **C** (Ф/м) и **G** (См/м) или, кратко, **RLCG**-параметрами. Вычисленные матрицы затем используются при решении телеграфных уравнений Хевисайда или производных от них для анализа целостности сигналов, получения временного отклика и других параметров [4]. При этом для минимизации затрат времени на получение всех первичных параметров МПЛП особый интерес представляет вычисление матрицы **C** (ёмкостной матрицы [5]), поскольку остальные три матрицы являются её производными [6].

В работе [7] представлен подробный вывод моделей вычисления ёмкостной матрицы с использованием MoM в сочетании с выражениями в замкнутом виде для элементов матрицы СЛАУ, для двумерных и трёхмерных структур с границами произвольной сложности, включающих идеально проводящую плоскость и без неё. В качестве операторного уравнения выступает уравнение:

$$\varphi = L^{-1}\rho, \ L^{-1} = \frac{1}{\varepsilon_0}\int G(\mathbf{r},\mathbf{r}')d\Gamma,$$

где $G(\mathbf{r}, \mathbf{r}')$ – соответствующая функция Грина, $\mathbf{r}$ – точка наблюдения $(x, y)$, $\mathbf{r}'$ – точка источника $(x', y')$, а $d\Gamma$ – дифференциал по поверхности структуры. В данной постановке задачи считаются заданными граничные условия по приложенному напряжению $\varphi$, при этом требуется найти плотность заряда $\rho$. Для двумерного случая функция Грина имеет вид:

$$G(\mathbf{r},\mathbf{r}') = -\frac{\ln|\mathbf{r}-\mathbf{r}'|}{2\pi},$$

а

$$\nabla G(\mathbf{r},\mathbf{r}') = -\frac{\mathbf{r}-\mathbf{r}'}{2\pi|\mathbf{r}-\mathbf{r}'|^2}.$$

Далее кратко рассмотрим математическую модель для вычисления ёмкостной матрицы на примере связанной микрополосковой линии (МПЛ), поперечное сечение которой приведено на рисунке 1. Структура содержит два проводника (I и II), которые расположены на диэлектрическом основании с относительной диэлектрической проницаемостью $\varepsilon_{r2}$ над идеально проводящей (бесконечной) плоскостью.

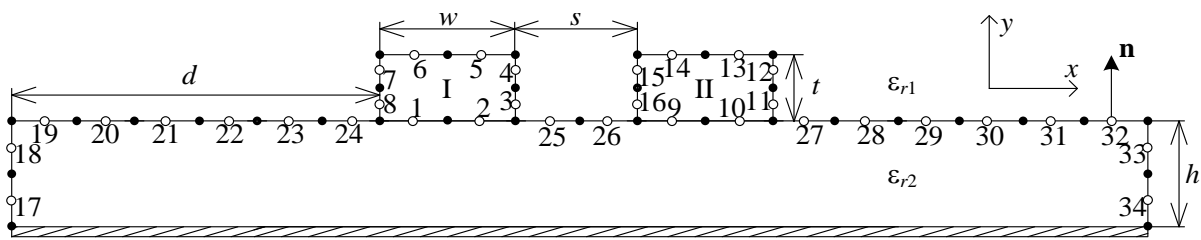

Рисунок 1. Поперечное сечение связанной МПЛ

Для вычисления необходимо пошагово выполнить следующие действия:

1. Дискретизировать границы проводников и диэлектрической подложки (граница раздела двух сред) на небольшие отрезки (подынтервалы) длинами $l_n$.

2. Когда в структуре имеются границы соприкосновения диэлектрика с проводником, необходимо работать в терминах полной плотности заряда $\sigma_T$, которая представляет собой сумму плотности свободного и поляризационного зарядов [8]: $\sigma_T(\mathbf{r}) = \sigma_S(\mathbf{r}) + \sigma_P(\mathbf{r})$. При этом, на границе диэлектрик-диэлектрик, полная плотность заряда состоит только из плотности поляризационных зарядов.

3. Учесть наличие плоскости земли с помощью метода зеркальных изображений.

4. Задать потенциалы (1 В) на проводниковых подынтервалах и получить результирующие уравнения для границ проводник-диэлектрик и диэлектрик-диэлектрик.

6. Аппроксимировать полную плотность заряда линейной комбинации из $N$ известных базисных функций.

7. Сформировать СЛАУ вида $\mathbf{S}\Sigma = \mathbf{V}$, где $\Sigma$ и $\mathbf{V}$ – матрицы размера $N \times N_{\text{COND}}$, где $\mathbf{S}$ – матрица размера $N \times N$.

8. Решить СЛАУ.

9. Вычислить элементы ёмкостной матрицы $\mathbf{C}$

$$C_{ij} = \int_{L_{Ci}} \varepsilon_r(\mathbf{r})\sigma_S^j(\mathbf{r})dl/V, V=1$$

где индекс $i$ – относится к проводнику, по контуру $L_{Ci}$ которого ведется интегрирование, а $j$ – к проводнику, находящемуся под потенциалом 1 В, когда остальные проводники под потенциалом 0 В. В матричном виде получим

$$C_{ij} = \sum_{k \in L_{Ci}} \varepsilon_r \Sigma_{kj} l_k /V, \ i,j = 1,...,N_{\text{COND}}, V=1\text{В}.$$

**2. Многовариантный анализ МПЛП**

Случайные вариации геометрических параметров МПЛП, обусловленные технологическим процессом, а также поиск оптимальных значений этих параметров делает необходимым многократно вычислять матрицы погонных параметров для достижения требуемого технического результата конечного изделия. Тогда, при многовариантном анализе с использованием приведенной в разделе 1 математической модели, требуется повторять п. 6–9. На практике число вариаций параметров МПЛП часто достаточно велико, что ведёт к существенному росту суммарных вычислительных затрат и тем самым затрудняет эффективное проектирование. Поэтому проектирование МПЛП посредством многовариантного анализа является нетривиальной задачей. Для решения этой проблемы возможны следующие подходы: переход к массивно-параллельным вычислениям, адаптация сетки к особенностям решения, разработка новых и совершенствование известных методов решения СЛАУ [9], их комбинация.

## 3. Методы построения сеток

По структуре получаемых сеток методы их построения можно условно разделить на методы построения равномерных и неравномерных сеток. Использование последовательного итерационного уточнения этих сеток приводит к, так называемым, адаптивным сеткам, которые, в большей степени, относятся к неравномерным сеткам [10]. При этом считается, что использование при расчетах таких сеток позволяет повысить их точность, не увеличивая числа узлов, тем самым экономя машинные ресурсы [11]. Поэтому рекомендуется сгущать сетку в подобластях резкого изменения контролируемой величины и тем самым повышать её точность, а для экономии машинных ресурсов использовать разреженную сетку в областях плавного изменения этой величины. Так, использование излишне густой сетки приводит к существенному увеличению затрат времени на вычисления. В связи с этим необходим оптимальный выбор сетки для получения требуемого решения с приемлемой точностью. Отдельно отметим, что существующие подходы в большей степени применимы к трехмерным задачам с треугольными или тетраэдральными сеткам и слабо применимы для двумерных задач и поверхностных подынтервальных сеток, как на рисунке 1 [12].

Применительно к методу моментов известны подходы к локальному учащению сетки с использованием нескольких правил оценки её качества [13, 14]. При этом, выбор их оптимальных параметров, сильно варьирующихся для разных структур, как и в случае с ручным учащением делает их слабо применимыми на практике при анализе МПЛП.

Известен подход, названный адаптивный итерационный выбор оптимальной сегментации (АИВОС), являющийся, по сути, представителем семейства методов *h*-уточнения, но основанный на равномерном (эквидистантном) учащении сетки. Его эффективность показана при анализе МПЛП [15]. Так, существенно сокращены затраты времени при незначительной потере точности. Сравнение выполнено с «ручной» равномерной сегментацией с длиной подынтервала равной трети от толщины проводника *t*/*n* (*n*=3 подынтервала на торце проводника). Из-за особенностей распределения зарядов по поверхности проводников и используемых кусочно-постоянных базисных функций такая длина подынтервалов часто считается допустимой для получения приемлемой точности решения [16]. Дальнейшее уменьшение длины подынтервалов ведет к росту вычислительных затрат на получение требуемых значений, при этом их уточнение часто оказывается малозначительным, что будет показано в следующем разделе.

Вернемся к рассмотрению методов из работ [13] и [15] (далее метод 1 и метод 2 соответственно). Далее приведем их общий псевдокод:

*Задать исходные параметры структуры (в т.ч. требуемую точность tol)*

*Задать начальную максимальную длину подынтервалов l*

*Задать максимальное допустимое число итераций $N_{it}^{max}$*

*Выполнить вычисление значения контролируемой величины $C_0$*

*Для i = 1, …, $N_{it}^{max}$*

   *Уменьшить длину подынтервалов*

   *Вычислить значение контролируемой величины $C_i$*

   *Если |$C_i$ – $C_{i–1}$| / |$C_{i–1}$| > tol*

      *Продолжить вычисление*

   *Иначе*

      *Выход из цикла*

*Увеличить i*

Прокомментируем стоку 6 псевдокода. При реализации метода 2 все подынтервалы, полученные на предыдущей итерации, разбиваются пополам. В

методе 1 подынтервалы также делятся пополам, но не все. Так, типовое значение составляет *p*=25 % от общего количества подынтервалов с максимальными значениями плотности заряда на них [13]. Главное различие в методах заключается в том, что первый даёт изначально неравномерную сетку, а второй – стремится к равномерной. Очевидно, что при стопроцентном учащении в методе 1, он становится эквивалентен методу 2.

В качестве контролируемой величины может выступать один из коэффициентов электростатической индукции, т.е. диагональный элемент матрицы **C** – $C_{kk}$, где *k* – порядковый номер проводника в структуре, как это сделано в [15], или норма Фробениуса матрицы **C** – [17].

Работа алгоритма начинается с использования грубой сетки, что требует минимальных затрат и является хорошим начальным приближением. Далее, итерационно происходит учащение сетки и анализ полученных результатов. Достоинством такого построения сетки является то, что оно позволяет в режиме реального времени контролировать ход моделирования. Если результаты от итерации к итерации изменяются несущественно (это контролируется параметром *tol*), то итерационный процесс останавливается. Также процесс остановится при достижении максимального допустимого числа итераций $N_{it}^{\max}$. Это необходимо для того, чтобы в случае отсутствия сходимости алгоритм корректно и за приемлемое время завершил свою работу.

Наконец, выбор начального шага сетки (длины подынтервалов *l*) оказывает существенное влияние на последующие суммарные затраты времени и памяти. Так, при слишком малом начальном шаге сетки эти затраты могут быть высоки, особенно при малой величине требуемой погрешности.

**4. Вычислительный эксперимент**

Для акцентирования внимания только на особенностях построения сетки для решения СЛАУ использован метод Гаусса. При вычислениях использованы система GNU Octave и персональный компьютер со следующими характеристиками: процессор – AMD Ryzen 3 3200G, тактовая частота – 3,6 ГГц; объём ОЗУ – 16 Гб; число ядер – 4.

Сначала оценено влияние равномерной сегментации вида *t*/*n* для трёх МПЛП (далее МПЛП 1 [15], МПЛП 2 и МПЛП 3 [18]), поперечные сечения которых приведены на рисунке 2. Параметры МПЛП 1 приведены в [15]. Параметры МПЛП 2: $t$=0,005 мм; $w$=0,05 мм; $s$=0,05 мм; $d$=0,15 мм; $h_1$=$h_2$=$h_3$=0,05 мм; $\varepsilon_{r1}$=$\varepsilon_{r3}$=3,8, $\varepsilon_{r2}$=2. Параметры МПЛП 3 приведены в работе [18]. В качестве контролируемых величин выступали собственные коэффициенты электростатической индукции контрольных проводников, обозначенных на рисунке 2 черным цветом (при нумерации этот проводник для всех МПЛП пронумерован первым), т. е. элементы $C_{11}$. Для сравнения эффективности методов использована сегментация вида *t*/3 являющаяся приемлемой для получения различия в точности контролируемой величины (элемент $C_{11}$) менее 1% [19].

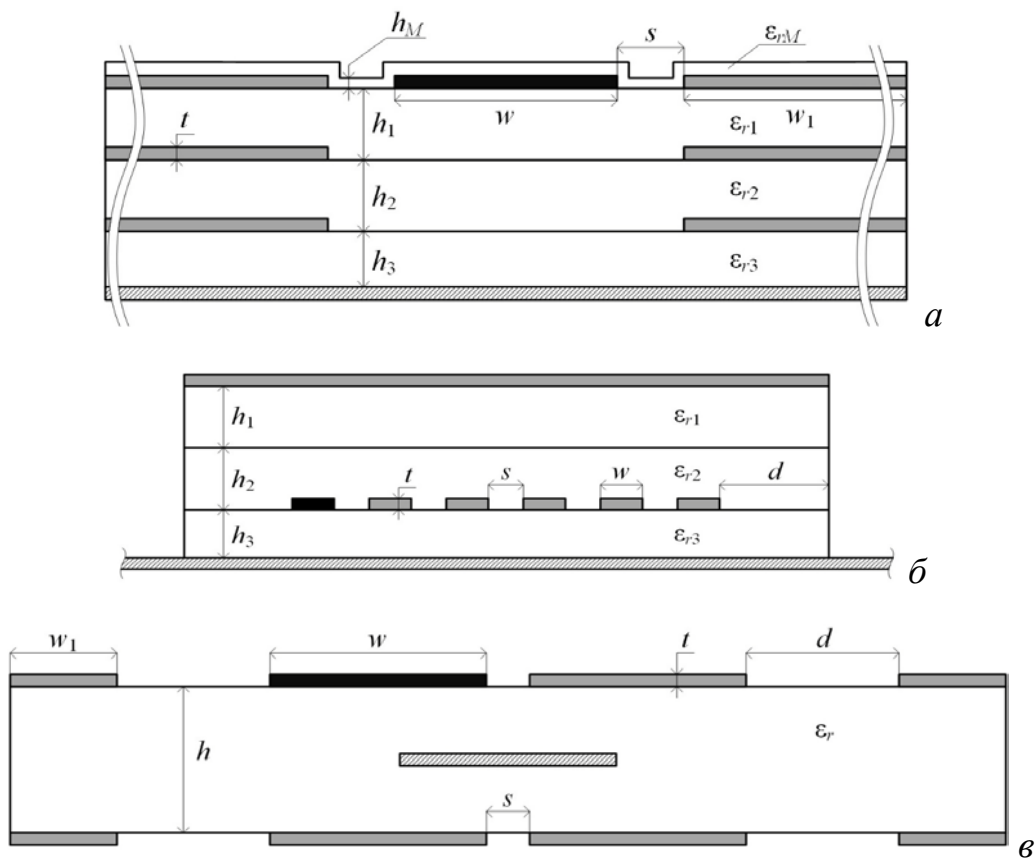

Рисунок 2. Поперечное сечение исследуемых структур: 1 (*а*), 2 (*б*) и 3 (*в*)

(серым и черным обозначены проводники, белым – диэлектрики)

Для МПЛП 1–3 произведены вычисления методами 1 и 2 при разной начальной длине подынтервалов, задаваемой как *l*=*kw*, где *k*=0,5; 1; 1,5; 2; 3.

Для метода 1 также изменялся процент разбиваемых подынтервалов $p=5; 15; 25\,\%$ (для метода 2 $p=100\,\%$). На основании предварительных вычислений было установлено, что оптимальными по точности и вычислительным затратам значениями параметра *tol* для методов 1 и 2 являются $10^{-3}$ и $10^{-2}$ соответственно. При вычислениях полагалось $N_{it}^{\max}=30$. В таблицу 1 сведены результаты вычислений: $\Delta C=|C-C_{t/3}|/|C_{t/3}|$ (в процентах), где $C_{t/3}$ и $C$ – значения $C_{11}$, полученные соответственно при сегментации $t/3$ и методом 1 или 2; отношения объема памяти $V_{t/3}$ и времени $T_{t/3}$ при сегментации $t/3$ к аналогичным значениям $V$ и $T$, полученным соответствующим методом 1 или 2 за число итераций $N_{it}$.

Из таблицы 1 видно, что метод 1 дает существенную экономию как машинной памяти, так и времени вычисления. При этом, чем меньше значение p, тем больше эта экономия. Однако, этот метод обладает для исследованных МПЛП достаточно низкой точностью ($\Delta C$ в большинстве случаев больше $1\,\%$ и доходит до $18\,\%$), что делает его не эффективным. Следует отметить, что корреляция параметров *l* с $\Delta C$ для данного метода мала и не линейна. Метод 2, напротив, показывает высокую зависимость от параметра *l*. Так, одновременно высокой точности ($\Delta C \leq 1\,\%$) и экономии машинных ресурсов (памяти и времени от 4 до 15 раз) удается достигнуть лишь при $k=2$. В остальных случаях либо слишком велика погрешность вычисления, либо отсутствует экономия машинных ресурсов. Поэтому далее использован метод 2 и $k=2$.

Для МПЛП 3 выполнена оценка эффективности метода 2 относительно сегментации при $t/3$ на примере раздельного изменения номинальных величин *t* и *w* в диапазоне *m* от минус $5\,\%$ до плюс $5\,\%$ (дисперсия параметров линии передачи [20]). Первые вычисления выполнены при $m=0\,\%$ (приведены в таблице 1), вторые – минус $5\,\%$, третьи – минус $4\,\%$ и т. д. (итого 11 вычислений). Полученные результаты (аналогично таблице 1) сведены в таблицу 2 кроме $N_{it}$, т. к. для всех значений *t* и *w* потребовалось 8 итераций.

Таблица 1. Результаты использования методов 1 и 2 относительно сегментации при $t/3$ для МПЛП 1–3

| $k$ | Метод | $p$, % | МПЛП 1 | | | | МПЛП 2 | | | | МПЛП 3 | | | |
|---|---|---|---|---|---|---|---|---|---|---|---|---|---|---|
| | | | $N_{it}$ | $\Delta C$ | $V_{t/3}/V$ | $T_{t/3}/T$ | $N_{it}$ | $\Delta C$ | $V_{t/3}/V$ | $T_{t/3}/T$ | $N_{it}$ | $\Delta C$ | $V_{t/3}/V$ | $T_{t/3}/T$ |
| 0,5 | 1 | 5 | 13 | 7,14 | 271,3 | 64,1 | 5 | 6,90 | 137,0 | 45,5 | 23 | 2,90 | 197,5 | 35,0 |
| | | 15 | 11 | 8,29 | 52,0 | 33,6 | 7 | 0,37 | 34,2 | 17,3 | 20 | 1,57 | 18,8 | 6,5 |
| | | 25 | 15 | 2,78 | 1,9 | 1,0 | 13 | 1,33 | 0,7 | 0,3 | 19 | 14,80 | 1,9 | 0,9 |
| | 2 | 100 | 6 | 0,18 | 1,9 | 1,9 | 2 | 19,36 | 62,4 | 65,6 | 9 | 0,51 | 0,7 | 0,5 |
| 1 | 1 | 5 | 19 | 1,07 | 224,5 | 76,6 | 10 | 0,07 | 204,9 | 40,5 | 18 | 17,92 | 453,1 | 107,5 |
| | | 15 | 17 | 1,44 | 28,2 | 14,6 | 12 | 0,43 | 25,0 | 10,4 | 21 | 1,59 | 15,4 | 5,3 |
| | | 25 | 17 | 2,57 | 2,2 | 1,2 | 9 | 0,45 | 15,2 | 9,3 | 19 | 14,80 | 1,9 | 0,9 |
| | 2 | 100 | 6 | 0,10 | 10,5 | 15,6 | 5 | 1,25 | 3,7 | 3,3 | 9 | 0,44 | 1,0 | 0,8 |
| 1,5 | 1 | 5 | 17 | 1,50 | 407,8 | 92,8 | 11 | 0,78 | 266,8 | 66,7 | 18 | 17,92 | 453,1 | 107,5 |
| | | 15 | 16 | 2,20 | 38,8 | 20,1 | 9 | 0,75 | 86,1 | 35,4 | 21 | 1,59 | 15,4 | 5,3 |
| | | 25 | 17 | 1,63 | 3,5 | 2,0 | 16 | 0,82 | 1,0 | 0,4 | 19 | 14,80 | 1,9 | 0,9 |
| | 2 | 100 | 6 | 0,15 | 24,1 | 26,7 | 5 | 1,71 | 7,7 | 7,5 | 8 | 0,36 | 6,7 | 6,0 |
| 2 | 1 | 5 | 20 | 6,15 | 315,4 | 80,7 | 10 | 0,99 | 403,8 | 47,4 | 18 | 17,92 | 453,1 | 107,5 |
| | | 15 | 27 | 3,25 | 3,0 | 1,3 | 10 | 0,61 | 106,5 | 30,6 | 21 | 1,59 | 15,4 | 5,3 |
| | | 25 | 19 | 0,89 | 1,5 | 0,7 | 7 | 0,84 | 107,3 | 35,2 | 19 | 14,80 | 1,9 | 0,9 |
| | 2 | 100 | 7 | 0,39 | 15,2 | 15,1 | 6 | 0,82 | 3,6 | 3,5 | 8 | 0,36 | 6,7 | 6,3 |
| 3 | 1 | 5 | 27 | 1,29 | 159,0 | 36,4 | 12 | 3,76 | 397,8 | 104,3 | 18 | 17,92 | 453,1 | 107,5 |
| | | 15 | 15 | 0,58 | 100,2 | 53,7 | 11 | 4,31 | 110,3 | 29,6 | 21 | 1,59 | 15,4 | 5,3 |
| | | 25 | 18 | 2,41 | 2,7 | 1,5 | 13 | 3,59 | 9,6 | 4,7 | 19 | 14,80 | 1,9 | 0,9 |
| | 2 | 100 | 7 | 0,20 | 24,1 | 24,6 | 6 | 1,43 | 9,3 | 9,2 | 2 | 62,27 | 18292,6 | 2784,8 |

Из таблицы 2 видно, что использование выявленных оптимальных, по критерию минимальных вычислительных затрат, параметров метода 1 позволяет существенно повысить производительность. Так, для 11 вычислений суммарный выигрыш по экономии машинной памяти составил до 6,7 раз, а затраченного времени – 8,3 раз относительно густой и равномерной сетки как при раздельном изменении $t$, так и – $w$. При этом очевидно, что в случае анализа при одновременном изменении этих параметров, экономия вычислительных затрат будет еще больше.

Таблица 2. Результаты использования метода 2 и $k=2$ относительно сегментации при $t/3$ для МПЛП 3 при изменении $t$ и $w$ в диапазонах ±5%

| $m$ | Изменение $t$ | | | | Изменение $w$ | | | |
|---|---|---|---|---|---|---|---|---|
| | $\Delta C$ | $N_Э/N$ | $V_Э/V$ | $T_Э/T$ | $\Delta C$ | $N_Э/N$ | $V_Э/V$ | $T_Э/T$ |
| –5 | 0,37 | 2,72 | 7,4 | 9,2 | 0,28 | 2,55 | 6,5 | 8,0 |
| –4 | 0,07 | 2,69 | 7,2 | 8,9 | 0,30 | 2,56 | 6,5 | 8,1 |
| –3 | 0,36 | 2,66 | 7,1 | 8,8 | 0,31 | 2,56 | 6,6 | 8,1 |
| –2 | 0,07 | 2,64 | 7,0 | 8,6 | 0,33 | 2,57 | 6,6 | 8,2 |
| –1 | 0,07 | 2,61 | 6,8 | 8,4 | 0,35 | 2,57 | 6,6 | 8,2 |
| 1 | 0,08 | 2,56 | 6,6 | 8,1 | 0,38 | 2,59 | 6,7 | 8,3 |
| 2 | 0,36 | 2,54 | 6,4 | 8,0 | 0,40 | 2,69 | 6,7 | 8,3 |
| 3 | 0,36 | 2,51 | 6,3 | 7,8 | 0,42 | 2,60 | 6,8 | 8,4 |
| 4 | 0,09 | 2,48 | 6,2 | 7,6 | 0,43 | 2,61 | 6,8 | 8,4 |
| 5 | 0,36 | 2,46 | 6,1 | 7,5 | 0,45 | 2,61 | 6,8 | 8,4 |

**Заключение**

Таким образом, в работе приведено наглядное описание математической модели для вычисления ёмкостной матрицы МПЛП, основанной на использовании квазистатического подхода и метода моментов. Кратко рассмотрены методы адаптивного построения сетки. Выполнено сравнение двух методов при изменении их параметров и выявлены их оптимальные значения по критерию минимальных вычислительных затрат относительно густой равномерной сетки. Это позволяет экономить вычислительные ресурсы при одновариантном и многовариантном анализе МПЛП в диапазоне её параметров. Так, на примере оценки влияния дисперсии параметров МПЛП на разброс значений её ёмкостной матрицы показано, что эта экономия существенна. При этом использованный метод достаточно просто интегрировать в существующие САПР. Далее целесообразно оценить его эффективность при оптимизации параметров МПЛП.

## Литература